\input amstex
\documentstyle{amsppt}
\topmatter \magnification=\magstep1 \pagewidth{5.2 in}
\pageheight{6.7 in}
\abovedisplayskip=10pt \belowdisplayskip=10pt
\parskip=8pt
\parindent=5mm
\baselineskip=2pt
\title
 On a $p$-adic interpolation function for the $q$-extension of the generalized Bernoulli polynomials
   and its derivative
\endtitle
\author   Taekyun Kim    \endauthor

\affil{ {\it Institute of Science Education,\\
        Kongju National University, Kongju 314-701, S. Korea\\
        e-mail: tkim64$\@$hanmail.net ( or tkim$\@$kongju.ac.kr)}}\endaffil
\thanks  2000 Mathematics Subject Classification:  11S80, 11B68, 11M99 .\endthanks
\abstract{ In [18], the new  $q$-extension of Bernoulli
polynomials and generalized Bernoulli numbers attached to $\chi$
were constructed by using $p$-adic invariant integral on $\Bbb
Z_p$. In this paper we construct the new $q$-extension of
generalized Bernoulli polynomials attached to $\chi$ due to author
and derive the existence of a specific $p$-adic interpolation
function which interpolate the $q$-extension of generalized
Bernoulli polynomials at negative integer. Finally, we give the
values of partial derivative  for this function and investigate
some properties which are related to this interpolation functions.
 }\endabstract
\rightheadtext{  On the $p$-adic interpolation function}
\leftheadtext{T. Kim}
\endtopmatter

\document

\head \S 1. Introduction \endhead Let $p$ be a fixed prime.
Throughout this paper  $\Bbb Z_p,\,\Bbb Q_p , \,\Bbb C$ and $\Bbb
C_p$ will, respectively, denote the ring of $p$-adic rational
integers, the field of $p$-adic rational numbers, the complex
number field and the completion of algebraic closure of $\Bbb Q_p
, $ cf.[3, 10, 17].  Let $v_p$ be the normalized exponential
valuation of $\Bbb C_p$ with $|p|_p=p^{-v_p(p)}=p^{-1}.$ When one
talks of $q$-extension, $q$ is variously considered as an
indeterminate, a complex number $q\in\Bbb C ,$ or a $p$-adic
number $q\in\Bbb C_p$. If $q\in\Bbb C ,$ one normally assumes
$|q|<1 .$  If $ q \in \Bbb C_p ,$  then we assume $|q-1|_p <
p^{-\frac1{p-1}},$ so that $q^x=\exp(x\log q)$ for $|x|_p \leq 1.$
Kubota and Leopoldt proved the existence of meromorphic functions,
$L_p(s, \chi )$, defined over the $p$-adic number field, that
serve as $p$-adic equivalents of the Dirichlet $L$-series.
 These $p$-adic
$L$-functions interpolate the values
$$L_p(1-n, \chi)=-\frac{1}{n}(1-\chi_n(p)p^{n-1})B_{n,\chi_n},
\text{ for $n\in\Bbb N=\{1, 2,\cdots, \}  $ ,}$$ where
$B_{n,\chi}$ denote the $n$th generalized Bernoulli numbers
associated with the primitive Dirichlet character $\chi ,$ and
$\chi_n=\chi w^{-n} ,$ with $w$  the $Teichm\ddot{u}ller$
character, cf.[1-24]. In a recent paper [18], the author
constructed the new $q$-extensions of Riemann zeta function,
Hurwitz's zeta function and Dirichlet's $L$-functions. In Section
2, we  define the $q$-extension of generalized Bernoulli
polynomials attached to $\chi$ and construct a new $q$-extension
of Hurwitz's type L-function which interpolates the $q$-extension
of generalized Bernoulli numbers attached to $\chi$ at negative
integers. In $\Bbb C$,  we introduce some of the basic facts about
the $q$-extension of Hurwitz's type $L$-function which
interpolates the generalized $q$-Bernoulli polynomials attached to
$\chi$. The values of this function at negative integers are
algebraic, hence may be regarded as lying in an extension of $\Bbb
Q_p$. We therefore look for a $p$-adic function which agrees with
at negative integers. The purpose of this paper is to construct
the new $q$-extension of generalized Bernoulli polynomials
attached to $\chi$ due to author and derive the existence of a
specific $p$-adic interpolation function which interpolate the
$q$-extension of generalized Bernoulli polynomials at negative
integer. Finally, we give the values of partial derivative  for
this $p$-adic function and investigate some properties which are
related to this $p$-adic interpolation functions in Section 3.

 \head 2. An $q$-extension of Dirichlet's $L$-function  \endhead

In this section we assume that $q\in\Bbb C$ with $|q|<1 $ and
$h\in\Bbb Z.$ The classical Bernoulli polynomials are  defined  by
$$\frac{t}{e^t-1}e^{xt}=\sum_{n=0}^{\infty}B_n(x)\frac{t^n}{n!},
\text{ $|t|<2\pi$ , }$$ and the classical Bernoulli numbers are
defined by $B_n=B_{n}(0).$ Let $\chi$ be a primitive Dirichlet
character  of conductor $f\in\Bbb N$. Then the generalized
Bernoulli polynomials associated with $\chi$, $B_{n,\chi}(x),$ are
also defined by
$$\sum_{a=1}^f\frac{\chi(a)te^{(a+x)t}}{e^{ft}-1}=\sum_{n=0}^{\infty}\frac{B_{n,\chi}(x)}{n!}t^n,
\text{ cf.[4, 18] }.\tag 1$$  Let us define the $q$-extension of
Bernoulli polynomials as follows:
$$\frac{h\log q +t
}{q^he^t-1}e^{xt}=\sum_{n=0}^{\infty}B_{n,q}^{(h)}(x)\frac{t^n}{n!},
\text{ cf.[4, 18] } . \tag2$$  We now consider the $q$-extension
of the generalized Bernoulli polynomials attached to $\chi$ as
$$\frac{\sum_{i=0}^{f-1}\left(te^{it}\chi(i)q^{hi}+e^{ti}\log
q^hq^{hi}\chi(i)\right)}{q^{hf}e^{ft}-1}e^{xt}=\sum_{n=0}^{\infty}B_{n,q,\chi}^{(h)}(x)\frac{t^n}{n!},
\text{ $|t|<2\pi $ }.\tag3$$ Let us define the $q$-extensions of
Bernoulli numbers and generalized Bernoulli numbers attached to
$\chi$ as $B_{n,q}^{(h)}=B_{n,q}^{(h)}(0)$ and
$B_{n,q,\chi}^{(h)}=B_{n,q}^{(h)}(0)$. Then we note that
$$B_{n,q,\chi}^{(h)}(x)=\sum_{k=0}^{n}\binom nk
B_{k,q,\chi}^{(h)}x^{n-k}=f^{n-1}\sum_{i=0}^{f-1}\chi(i)q^{hi}B_{n,q^f}^{(h)}(\frac{i+x}{f}).\tag4$$
Remark. Let $\delta_{1,k}$ be denoted by Kronecker symbol. Then we
see that $$B_{0,q}^{(h)}=\frac{h\log q}{q^h-1},\text{ }
q^h(B_q^{(h)}+1)^n-B_{n,q}^{(h)}=\delta_{1,k}, \text{ $k\geq 1$,
cf.[1, 2, 21, 22,]  },$$ with the usual convent of replacing
$(B_q^{(h)})^n$ by $B_{n,q}^{(h)}.$

For $s\in\Bbb C$, the $q$-extension of Hurwitz zeta function is
defined by
$$\zeta_{q}^{(h)}(s,x)=\sum_{n=0}^{\infty}\frac{q^{nh}}{(n+x)^s}-\frac{h\log
q}{s-1}\sum_{n=0}^{\infty}\frac{q^{nh}}{(n+x)^{s-1}}, \text{
cf.[12, 14, 18].}$$ Note that $\zeta_q^{(h)}(s,x)$ is analytic
continuation for $\Re(s)>1 .$ This function have the below value
at negative integers:
$$\zeta_q^{(h)}(1-n,x)=-\frac{B_{n,q}^{(h)}(x)}{n}, \text{ for
$n\in\Bbb N=\{1,2,\cdots \}$ }.$$ In [18], the $q$-extension of
Dirichlet's $L$-function is also defined by
$$L_q^{(h)}(s,\chi)=\sum_{n=1}^{\infty}\frac{\chi(n)q^{hn}}{n^s}-\frac{h\log
q}{s-1}\sum_{n=1}^{\infty}\frac{q^{nh}\chi(n)}{n^{s-1}}. \tag5$$
This function is analytic continuation  for $\Re(s)>1$. Note that
$L_q^{(h)}(1-n,\chi)=-\frac{B_{n,q,\chi}^{(h)}}{n}, $ for
$n\in\Bbb N,$ cf.[12, 14, 18, 19]. We now set
$$F_{q,\chi}^{(h)}(t,x)=\frac{\sum_{i=0}^{f-1}\left(te^{it}\chi(i)q^{hi}+e^{ti}\log
q^hq^{hi}\chi(i)\right)}{q^{hf}e^{ft}-1}e^{xt}, \text{ for
$|t|<\frac{2\pi}{f}$} .\tag6$$ By using (6), we easily see that
$$F_{q,\chi}^{(h)}(t,x)=-t\sum_{n=0}^{\infty}\chi(n)q^{nh}e^{(n+x)t}-h\log
q\sum_{n=0}^{\infty}\chi(n)q^{hn}e^{(n+x)t}.\tag7$$ It is easy to
see that the series on the right-hand side of (7) are uniformly
convergent. Hence, we have
$$
B_{k,q,\chi}^{(h)}(x)=\frac{d^k}{dt^k}F_{q,\chi}^{(h)}(t,x)|_{t=0}=-k\sum_{n=0}^{\infty}\chi(n)q^{hn}(n+x)^{k-1}
-h\log q\sum_{n=0}^{\infty}\chi(n)q^{hn}(n+x)^k.\tag8$$ That is,
$$-\frac{B_{k,q,\chi}(x)}{k}=\sum_{n=0}^{\infty}\chi(n)q^{hn}(n+x)^{k-1}+\frac{h\log
q}{k}\sum_{n=0}^{\infty}\chi(n)q^{hn}(n+x)^k, \text{ $k\in\Bbb
N$.} \tag9$$ Thus, we can consider the $q$-extension of
Dirichlet's $L$-function which interpolates the generalizes
q-Bernoulli numbers at negative integer as follows: \proclaim{
Definition 1} For $s\in\Bbb C ,$ define
$$L_q^{(h)}(s,x|\chi)=\sum_{n=0}^{\infty}\frac{\chi(n)q^{hn}}{(n+x)^s}-\frac{h\log
q}{s-1}\sum_{n=0}^{\infty}\frac{q^{hn}\chi(n)}{(n+x)^{s-1}}.\tag10$$
\endproclaim
Note that $L_q^{(h)}(s,x|\chi)$ is analytic continuation in $\Bbb
C$ with only simple pole at $s=1.$

By (9) and (10), we obtain the following:
 \proclaim{ Proposition 2}
For any positive integer $k$, we have
$$L_q^{(h)}(1-k,x|\chi)=-\frac{B_{k,q,\chi}^{(h)}}{k}(x).$$
\endproclaim
We now give the integral representation of the $q$-extension of
Dirichlet's $L$-function which interpolates the $q$-extension of
generalized Bernoulli polynomials attached to $\chi$ in $\Bbb C$.
Let $\Gamma(s)$ be the gamma function. Then we can readily see
that
$$L_q^{(h)}(s,x|\chi)=\frac{1}{\Gamma(s)}\int_{0}^{\infty}t^{s-2}F_{q,\chi}^{(h)}(-t,x)dt. \tag11$$
From (3),(6) and (11), we note that
$$L_q^{(h)}(1-n, x|\chi)=-\frac{B_{n,q,\chi}^{(h)}(x)}{n}, \text{  for
$n\in\Bbb N$.}$$ Let $s$ be a complex variable , $a$ and $F$ be
integers with $0<a<f$ Then we define $H_q^{(h)}(s,a|F)$ as
follows:
$$
H_q^{(h)}(s,a|f)=\sum_{\Sb m\equiv a(\mod f)\\m>0 \endSb}
\frac{q^{mh}}{m^s}-\frac{\log q}{s-1}\sum_{\Sb m\equiv a(\mod
f)\\m>0\endSb}\frac{q^{hm}}{m^{s-1}}=f^{-s}q^{ha}\zeta_{q^f}^{(h)}(s,\frac{a}{f}).
\tag12$$ Let $\chi(\neq 1)$ be the Dirichlet's character with
conductor $f\in\Bbb N$. Then the $q$-analogue of Dirichlet's
$L$-function can be expressed as the sum.
$$L_q^{(h)}(s,\chi)=\sum_{a=1}^{f}\chi(a)H_q^{(h)}(s,a|f), \text{
for $s\in\Bbb C$.}\tag13$$ The function $H_q^{(h)}(s,a|f)$ is a
meromorphic for $s\in\Bbb C$ with simple pole at $s=1$, having
residue $\frac{q^{ha}\log q}{q^{hf}-1}$, and it interpolates the
values
$$H_q^{(h)}(1-n,a|f)=-\frac{f^{n-1}}{n}q^{ha}B_{n,q^f}(\frac{a}{f}), \text{ where $n\in\Bbb Z,$ $n\geq 1$.}\tag14$$
Now, we modify the $q$-analogue of the partial zeta function  as
follows:
$$H_q^{(h)}(s,a|f)=\frac{q^{ha}}{(s-1)f}a^{1-s}\sum_{j=0}^{\infty}\binom{1-s}{j}\left(\frac{f}{a}\right)^j
B_{j,q^f}^{(h)}, \text{ for $s\in\Bbb C$ .} \tag15$$ From (13),
(14) and (15), we can derive the below:
$$L_q^{(h)}(s,\chi)=\frac{1}{s-1}\frac{1}{f}\sum_{a=1}^{f}\chi(a)q^{ha}a^{1-s}\sum_{m=0}^{\infty}
\binom{1-s}{m}\left(\frac{f}{a}\right)^mB_{m,q^f}^{(h)}.\tag16$$
It is easy to see that we can express $L_q^{(h)}(s,x|\chi)$ in the
similar method. Using (15) to define $H_q^{(h)}(s,a+x|f)$ for all
$a\in\Bbb Z$ with $0<a<f, $ $x\in\Bbb R$ with $0<x<1,$ we can
define
$$L_q^{(h)}(s,x|\chi)=\sum_{a=1}^f \chi(a)q^{-hx}H_q^{(h)}(s, a+x|f).$$
Let $f$ and $a$ be the positive integers with $0<a<f.$ Then we
have
$$L_q^{(h)}(s,x|\chi)=\frac{1}{s-1}\frac{1}{f}\sum_{a=1}^f\chi(a)q^{ha}(a+x)^{1-s}
\sum_{m=0}^{\infty}\binom{1-s}{m}\left(\frac{f}{a+x}\right)^m
.\tag17$$ By (17), we see that $L_q^{(h)}(s,x|\chi)$ is an
analytic for $x\in\Bbb R$ with $0<x\leq 1 ,$ $s\in\Bbb C$, except
$s\neq 1 .$ Furthermore, for each $n\in\Bbb Z$ with $n\geq 1$, we
have
$$L_q^{(h)}(1-n,x|\chi)=-\frac{B_{n,q,\chi}^{(h)}(x)}{n}.\tag18$$ In
$\Bbb C$,  we introduced some of the basic facts about the
$q$-extension of Dirichlet $L$-function which interpolates the
generalized $q$-Bernoulli polynomials attached to $\chi$. The
values of $L_q^{(h)}(s,x|\chi)$ at negative integers are
algebraic, hence may be regarded as lying in an extension of $\Bbb
Q_p$. We therefore look for a $p$-adic function which agrees with
at negative integers in the next section.

 \head 3. On a $p$-adic interpolation function for the $q$-extension of the generalized Bernoulli polynomials
   and its derivative  \endhead
In this section we shall consider the $p$-adic  analogs of the
$q$-$L$-functions, $L_q^{(h)}(s,x|\chi),$ which were introduced in
the previous section. Indeed this functions are the $q$-analogs of
the $p$-adic interpolation functions for the generalized Bernoulli
polynomials attached to $\chi$. Let $w$ denote the
$Teichm\ddot{u}ller$ character, having conductor $f_w =p^*$. For
an arbitrary character $\chi$, we define $\chi_n=\chi w^{-n},$
where $n\in\Bbb Z$, in the sense of the product of characters.
Throughout this section, we assume that $q\in\Bbb C_p$ with
$|1-q|_p <p^{-\frac{1}{p-1}} .$ Let
$<a>=w^{-1}(a)a=\frac{a}{w(a)}.$ Then, we note that $<a>\equiv 1$
$(\mod p^*\Bbb Z_p).$ By the definition of $<a>$, we easily see
that
$<a+p^*t>=w^{-1}(a+p^*t)(a+p^*t)=w^{-1}(a)a+w^{-1}(a)(p^*t)\equiv
1$ $(\mod p^*\Bbb Z_p[t]),$ where $t\in\Bbb C_p$ with $|t|_p\leq
1,$ $(a,p)=1 .$ The $p$-adic logarithm function, $\log_p ,$ is the
unique function $\Bbb C_p^{\times} \rightarrow \Bbb C_p $ that
satisfy (1) $\log_p (1+x)=\sum_{n=1}^{\infty}\frac{(-1)^n}{n}x^n
,$ $|x|_p<1$, (2) $\log_p(xy)=\log_p x+\log_p y ,$ $\forall x,
y\in\Bbb C_p^{\times},$ and $\log_p p =0 .$ Let
$A_j(x)=\sum_{n=0}^{\infty}a_{n,j}x^n$, $a_{n,j}\in\Bbb C$, $j=0,
1, 2,\cdots$ be a sequence of power series, each of which
converges in a fixed subset $D=\{s\in\Bbb C_p||s|_p\leq
|p^*|^{-1}p^{-\frac{1}{p-1}}\}$ of $\Bbb C_p$ such that (1)
$a_{n,j}\rightarrow a_{n, 0}$ as $j\rightarrow \infty$ for
$\forall n$; (2) for each $s\in D$ and $\epsilon >0$, there exists
$n_0=n_0(s,\epsilon)$ such that $\left|\sum_{n\geq
n_0}a_{n,j}s^n\right|_p<\epsilon $ for $\forall j$. Then
$\lim_{j\rightarrow \infty}A_j(s)=A_0(s)$ for all $s\in D .$ This
is used by Washington [24] to show that each of the function
$w^{-s}(a)a^s$ and $\sum_{m=0}^{\infty}\binom sm
\left(\frac{F}{a}\right)^mB_m ,$ where $F$ is the multiple of
$p^*$ and $f=f_{\chi}$, is analytic in $D$.
 Let $F$ be a positive
integral multiple of $p^*$ and $f=f_{\chi}$, and let
$$\aligned
&L_{p,q}^{(h)}(s,t|\chi) \\
&=\frac{1}{s-1}\frac{1}{F}\sum_{\Sb a=1\\
(a,p)=1\endSb}^F\chi(a)q^{ha}<a+p^*t>^{1-s}
 \sum_{m=0}^{\infty} \binom{1-s}m
\left(\frac{F}{a+p^*t}\right)^mB_{m,q^F}^{(h)}.\endaligned\tag19$$
 Then $L_{p,q}^{(h)}(s,t|\chi)$
is analytic for $t\in\Bbb C_p$ with $|t|_p\leq 1 ,$ provided $s\in
D$, except $s\neq 1$ when $\chi\neq 1.$
 For $t\in\Bbb C_p$ with
$|t|_p\leq 1$, we see that $\sum_{j=0}^{\infty}\binom sj
\left(\frac{F}{a+p^*t}\right)^jB_{j,q^F}^{(h)}$ is analytic for
$s\in D .$ It readily follows that
$<a+p^*t>^s=<a>^s\sum_{m=0}^{\infty}\binom sm
\left(a^{-1}p^*t\right)^m $ is analytic for $t\in\Bbb C_p$ with
$|t|_p\leq 1$ when $s\in D .$ Thus, since
$(s-1)L_{p,q}^{(h)}(s,t|\chi)$ is a finite sum of products of
these two functions, it must also be analytic for $t\in\Bbb C_p$,
$|t|_p\leq 1$, whenever $s\in D.$ Note that
$$\lim_{s\rightarrow 1}(s-1)L_{p,q}^{(h)}(s,t|\chi)=\frac{1}{F}\sum_{\Sb a=1\\(a,p)=1\endSb}^{F}
\chi(a)q^{ha}B_{0,q^F}^{(h)}. \tag20$$ We now let $n\in\Bbb Z,$ $
n\geq 1, $ and fix $t\in\Bbb C_p$ with $|t|_p\leq 1$. Since $F$
must be a multiple of $f=f_{\chi_n}$. By (4), we see that
$$B_{n,q,\chi_n }^{(h)}(p^*t)=F^{n-1}\sum_{a=0}^{F-1}\chi_n(a)q^{ha}B_{n,q^F}^{(h)}(\frac{a+p^*t}{F}).
\tag 21$$ If $\chi_n(p)=0$, then $(p, f_{\chi_n})=1$, so that
$\frac{F}{p}$ is a multiple of $f_{\chi_n}$. Therefore, we obtain
$$\chi_n(p)p^{n-1}B_{n,q^p,\chi_n}^{(h)}(p^{-1}p^*t)
=F^{n-1}\sum_{\Sb a=0\\p|a
\endSb}^F\chi_n(a)q^{ha}B_{n,q^F}^{(h)}(\frac{a+p^*t}{F}). \tag 22$$
The difference of these quantities  yields
$$B_{n,q,\chi_n}^{(h)}(p^*t)-\chi_n(p)p^{n-1}B_{n,q^p,\chi_n}^{(h)}(p^{-1}p^*t)
=F^{n-1}\sum_{\Sb a=1\\ p\nmid
a\endSb}^F\chi_n(a)q^{ha}B_{n,q^F}(\frac{a+p^*t}{F}).\tag23$$ From
the definition of $q$-Bernoulli polynomials (see[18:p.3]), we note
that
$$\aligned
B_{n,q^F}^{(h)}(\frac{a+p^*t}{F})&=\sum_{k=0}^{n}\binom nk
\left(\frac{a+p^*t}{F}\right)^{n-k}B_{k,q^F}^{(h)}\\
&=F^{-n}(a+p^*t)^n\sum_{k=0}^n \binom nk
\left(\frac{F}{a+p^*t}\right)^k B_{k,q^F}^{(h)}.
\endaligned$$
Since $\chi_n(a)=\chi(a)w^{-n}(a)$ and for $(a,p)=1$, and
$t\in\Bbb C_p $ with $|t|_p\leq 1$, we have
$$\aligned
&B_{n,q,\chi_n}^{(h)}(p^*t)-\chi_n(p)p^{n-1}B_{n,q^p,\chi_n}^{(h)}(p^{-1}p^*t)\\
&=\frac{1}{F}\sum_{\Sb a=1\\ (a,p)=1\endSb }^F\chi(a) q^{ha}
<a+p^*t>^n\sum_{m=0}^{\infty}\binom nm
\left(\frac{F}{a+p^*t}\right)^m B_{m,q^F}^{(h)}
\endaligned\tag24$$
Thus, we see that
$$-\frac{1}{n}\left(B_{n,q,\chi_n}^{(h)}(p^*t)-\chi_n(p)p^{n-1}B_{n,q^p, \chi_n}^{(h)}(p^{-1}p^*t)\right)
=L_{p,q}^{(h)}(1-n,t|\chi),\text{ for $ n\in\Bbb N.$} $$ Therefore
we obtain the following theorem: \proclaim{ Theorem 3} Let $F$ be
a positive integral multiple of $p^*$ and $f=f_{\chi}$, and let
$$\aligned
&L_{p,q}^{(h)}(s,t|\chi)\\
&=\frac{1}{s-1}\frac{1}{F}\sum_{\Sb
a=1\\(a,p)=1\endSb}^F\chi(a)q^{ha}<a+p^*t>^{1-s}\sum_{m=0}^{\infty}\binom{1-s}m
\left(\frac{F}{a+p^*t}\right)^m B_{m,q^F}^{(h)} .\endaligned $$
Then, $L_{p,q}^{(h)}(s, t|\chi)$ is analytic for $t\in\Bbb C_p$,
$|t|_p\leq 1,$ $h\in\Bbb Z$, provided $s\in D$, except $s\neq 1$.
Also, if $t\in\Bbb C_p ,$ $|t|_p\leq 1$, this function is analytic
for $s\in D$ when $\chi \neq 1$, and meromorphic for $s\in D$,
with simple pole at $s=1$ having residue $\frac{h\log_p
q}{q^{hF}-1}\left(\frac{1-q^{hF}}{1-q^h}-\frac{1-q^{hF}}{1-q^{hp}}\right)
$ when $\chi=1$. Furthermore, for each $n\in\Bbb Z$, $n\geq 1$, we
have
$$L_{p,q}^{(h)}(1-n,t|\chi)=-\frac{1}{n}\left(B_{n,q,\chi_n}^{(h)}(p^*t)-\chi_n(p)p^{n-1}
B_{n,q^p,\chi_n}^{(h)}(p^{-1}p^*t)\right).$$
\endproclaim
\proclaim {Remark} Note that $\lim_{h\rightarrow
0}L_{p,q}^{(h)}(s, 0|\chi)=\lim_{q\rightarrow
1}L_{p,q}^{(h)}(s,0|\chi)=L_p(s,\chi)$ for $s\in D$ with $s\neq 1$
if $\chi =1,$ where $L_{p}(s,\chi)$ is Kubota-Leopoldt's $p$-adic
$L$-function, cf.[5, 6, 7, 8, 11, 13, 20, 23, 24].
\endproclaim
We now consider the $q$-analogue of the partial $p$-adic zeta
function as follows:
$$H_{p,q}^{(h)}(s,a|F)=\frac{1}{s-1}\frac{1}{F}<a>^{1-s}\sum_{j=0}^{\infty}
\binom{1-s}j \left(\frac{F}{a}\right)^j B_{j,q^F}^{(h)},
$$ where $s\in D$, $s\neq 1$, $a\in\Bbb Z$ with $(a,p)=1$, and $F$
is a multiple of $p^*$, cf. [8, 13]. The function
$L_{p,q}^{(h)}(s,\chi)$ can be rewritten as the sum
$$L_{p,q}^{(h)}(s,\chi)=\sum_{\Sb a=1\\(a,p)=1
\endSb}^F\chi(a)q^{ha}H_{p,q}^{(h)}(s,a|F), \text{ cf.[ 8, 13], }$$
provided $F$ is a multiple of both $p^*$ and $f=f_{\chi} .$ The
function $H_{p,q}^{(h)}(s,a| F)$ is a meromorphic for $s\in D$
with a simple pole at $s=1$, having residue $\frac{h\log_p
q}{q^{hF}-1} ,$ and it interpolates the values
$$H_{p,q}^{(h)}(1-n,a|F)=-\frac{1}{n}F^{n-1}w^{-n}(a)B_{n,q^F}^{(h)}(\frac{a}{F}),
$$ where $n\in\Bbb Z$, $n\geq 1$, cf. [8-19].
By using $H_{p,q}^{(h)}(s, a+p^*t| F)$, we can express
$L_{p,q}^{(h)}(s,t|\chi)$ for all $a\in\Bbb Z$, $(a,p)=1$, and
$t\in\Bbb C_p$ with $|t|_p\leq 1$, as follows:
$$L_{p,q}^{(h)}(s,t|\chi)=\sum_{\Sb a=1\\(a,p)=1
\endSb}^F\chi(a)q^{ha}H_{p,q}^{(h)}(s,a+p^*t| F). $$
It is easy to see that $H_{p,q}^{(h)}(s, a+p^*t| F)$ is analytic
for $t\in\Bbb C_p$, $|t|_p\leq 1$, where $s\in D$, $s\neq 1$, and
meromorphic for $s\in D$, with a simple pole at $s=1$, when
$t\in\Bbb C_p$, $|t|_p\leq 1$. Let us consider the first partial
derivative of the function $L_{p,q}^{(h)}(s,t|\chi)$ at $s=0$.
 The value of
$\frac{\partial}{\partial s}L_{p,q}^{(h)}(0,t|\chi)$ is the
coefficient of $s$ in the expansion of $L_{p,q}^{(h)}(s,t|\chi)$
at $s=0 .$ By using Taylor expansion at $s=0$, we see that
$$\aligned
&\frac{1}{1-s}=1+s+\cdots, \\
&<a+p^*t>^{1-s}=<a+p^*t>\left(1-s\log_p <a+p^* t>+\cdots \right),\\
&\binom{1-s}{m}=\frac{(-1)^{m+1}}{m(m-1)}s +\cdots .
\endaligned$$
By employing these expansion, along with some algebraic
manipulation, we evaluate $\frac{\partial}{\partial
s}L_{p,q}^{(h)}(0, t|\chi).$ From the definition of
$L_{p,q}^{(h)}(s,t|\chi)$, we note that
$$\aligned
&L_{p,q}^{(h)}(s,t|\chi)\\
&=\frac{1}{s-1}\frac{1}{F}\sum_{\Sb a=1
\\(a, p)=1 \endSb}^F\chi(a)q^{ha}<a+p^*t>^{1-s}\sum_{m=0}^{\infty}
\binom{1-s}{m}\left(\frac{F}{a+p^*t}\right)^mB_{m,q^F}^{(h)}.
\endaligned$$
Thus, we have
$$\aligned
&\frac{\partial}{\partial s}L_{p,q}^{(h)}(s,t|\chi)|_{s=0}
=\sum_{a=1}^F\chi_1(a)q^{ha}\big\{\left(\left(\frac{a+p^*t}{F}\right)B_{0,q^F}^{(h)}+B_{1,q^F}^{(h)}\right)\log_p
<a+p^*t>\\
&-\left(\frac{a+p^*t}{F}\right)B_{0,q^F}^{(h)}
+\sum_{m=2}^{\infty}\frac{(-1)^m}{m(m-1)}
B_{m,q^F}^{(h)}\left(\frac{a+p^*t}{F}\right)^{1-m}\big\}
-\sum_{\Sb a=1\\(a,p)=1\endSb}^F\chi_1(a)q^{ha}B_{1,q^F}^{(h)}.
\endaligned\tag25$$
 Since the Diamond gamma function is defined by
$$G_p(x)=(x-\frac{1}{2})\log_p
x-x+\sum_{j=2}^{\infty}\frac{B_j}{j(j-1)}x^{1-j}, \text{ for
$|x|_p>1$, cf. [5, 6, 20, 23] },\tag26$$ and $w(a)$ is a root of
unity for $(a, p)=1$, we see that
$$ \log_p<a+p^*t>=\log_p(a+p^*t)+\log_p w^{-1}(a)=\log_p(a+p^*t), \text{ cf. [24]}.$$
For $f\in UD(\Bbb Z_p, \Bbb C_p)=\{f|f:\Bbb Z_p \rightarrow \Bbb
C_p \text{ is uniformly differentiable function }\}$, the
Volkenborn integral is defined by
$$I_0(f)=\int_{\Bbb Z_p}f(x)d\mu_0(x)=\lim_{N\rightarrow
\infty}\frac{1}{p^N}\sum_{x=0}^{p^N-1} f(x), \text{ cf.
[18].}\tag27$$ It is easy to see that
$$I_0(f_1)=I_0(f)+f^{\prime}(0), \text{ where $f_1(x)=f(x+1).$}$$
If we take $f(x)=q^{hx}e^{tx},$ ($h\in\Bbb Z$), then we have
$$\int_{\Bbb Z_p}q^{hx}e^{xt}d\mu_0(x)=\frac{h\log_p q
+t}{q^he^t-1}. $$ Thus, we note that
$$\int_{\Bbb Z_p}q^{hx}x^nd\mu_0(x)=B_{n,q}^{(h)}, \text{ for
$h\in\Bbb Z$ }.\tag 28$$
 We now consider a locally analytic function $G_{p,q}^{(h)}(x)$ which are the $q$-extension  of Diamond gamma
 function as follows:
 $$G_{p,q}^{(h)}(x)=\int_{\Bbb
 Z_p}\left\{(x+z)\log_p(x+z)-(x+z)\right\}q^{hz}d\mu_0(z),
 \text{ for $|x|_p>1. $}\tag29 $$
From the above Eq.(29), we note that $G_{p,q}^{(h)}(x)$ is locally
analytic on $\Bbb C_p\setminus\Bbb Z_p $ . By (28) and (29), we
easily see that
$$G_{p,q}^{(h)}(x)=\left(xB_{0,q}^{(h)}+B_{1,q}^{(h)}\right)\log_p
x-xB_{0,q}^{(h)}+\sum_{n=1}^{\infty}\frac{(-1)^{n+1}}{n(n+1)}B_{n+1,q}^{(h)}x^{-n},\text{
for $|x|_p>1 $.}\tag30$$ Note that $\lim_{q\rightarrow
1}G_{p,q}^{(h)}(x)=\lim_{h\rightarrow 0}G_{p,q}^{(h)}(x)=G_p(x).$
From (25) and (30), we can derive the below formula:
$$\aligned \frac{\partial}{\partial s}L_{p,q}^{(h)}(0,t|\chi)
&=\sum_{\Sb a=1
\\(a,p)=1\endSb}^F\chi_1(a)q^{ha}G_{p,q^F}^{(h)}\left(\frac{a+p^*t}{F}\right)
-L_{p,q}^{(h)}(0,\chi)\log_p F\\& -\sum_{\Sb a=1\\(a,p)=1
\endSb}^F\chi_1(a)q^{ha}B_{1,q^F}^{(h)}.
\endaligned$$
Therefore we obtain the following theorem:
 \proclaim{ Theorem 4}
 Let $\chi$ be the primitive Dirichlet character, and let $F$ be a
 positive integral multiple of $p^*$ and $f=f_{\chi}$. Then for
 any $t\in\Bbb C_p$ with $|t|_p\leq 1$, $h\in\Bbb Z$, we have
$$\aligned
&\frac{\partial}{\partial s}L_{p,q}^{(h)}(0,t|\chi) =\sum_{\Sb a=1
\\(a,p)=1\endSb}^F\chi_1(a)q^{ha}G_{p,q^F}^{(h)}(\frac{a+p^*t}{F})
-L_{p,q}^{(h)}(0,\chi)\log_p F \\& -\sum_{\Sb a=1\\(a,p)=1
\endSb}^F\chi_1(a)q^{ha}B_{1,q^F}^{(h)}.
\endaligned$$ \endproclaim
Note that $$\lim_{h\rightarrow 0}\frac{\partial}{\partial
s}L_{p,q}^{(h)}(0,t|\chi)= \frac{\partial}{\partial
s}L_{p}(0,t|\chi)=\sum_{\Sb a=1
\\(a,p)=1\endSb}^F\chi_1(a)G_{p}(\frac{a+p^*t}{F})
-L_{p}(0,\chi)\log_p F. $$ This formula can be considered  as the
generalization of Ferrero-Greenberg's theorem, cf.[5, 6, 7, 8, 11,
13, 20, 23, 24].

\enddocument